\documentclass{amsart}
\usepackage[english]{babel}
\usepackage{hyperref}
\usepackage{cleveref}
\usepackage{amssymb, amsfonts, amsmath}
\usepackage{pgfplots}
\usepackage{comment} 
\usepackage{mathtools}
\usepackage{wrapfig}
\usepackage{tikz}
\usepackage{comment}
\usepackage{mathrsfs}
\usepackage{tabularx, hyperref}
\usepackage{amssymb} \usepackage{amsfonts} \usepackage{amsmath}
\usepackage{amsthm} \usepackage{epsfig, subfig}
\usepackage{ amscd, amsxtra, latexsym}
\usepackage{epsfig,  graphicx, psfrag}
\usepackage[all]{xy}
\usepackage{caption}
\usepackage{enumerate}
\usepackage{color}
\DeclareMathOperator{\vol}{Vol}

\newcommand{\str}{{{\mathrm{str}}}}

\newcommand{\tilM}{\widetilde{M}}

\newcommand\restr[2]{{
\left.\kern-\nulldelimiterspace 
#1 
\vphantom{\big|} 
\right|_{#2} 
}}

\newcommand\rrestr[2]{{
\left.\kern-\nulldelimiterspace 
\left.\kern-\nulldelimiterspace 
#1 
\vphantom{\big|} 
\right|\hspace{-2.4pt} 
\right|_{#2} 
}}

\renewcommand{\phi}{\varphi}
\renewcommand{\colon }{\,:}

\newcommand*{\Simpl}{\mathrm{Sing}}
\newcommand*{\Comb}{\mathrm{Vert}}
\newcommand*{\Chain}{C^\mathrm{Vert}}
\newcommand*{\ChainSing}{C^\mathrm{Sing}}
\newcommand*{\Cochain}{C_{\mathrm{Vert},b}}
\newcommand*{\CochainSing}{C_{\mathrm{Sing},b}}
\newcommand*{\STR}{\mathrm{str}}

\newtheorem{lemma}{Lemma}[section]

\newtheorem*{prop*}{Proposition}
\newtheorem{prop_intro}{Proposition}

\newtheorem{thm_intro}[prop_intro]{Theorem}

\newtheorem{cor_intro}[prop_intro]{Corollary}

\theoremstyle{definition}

\newtheorem{defn}[lemma]{Definition}

\theoremstyle{remark}
\newtheorem{rmk}[lemma]{Remark}

\newtheoremstyle{citing}
{3pt}
{3pt}
{\itshape}
{}
{\bfseries}
{.}
{.5em}
{\thmnote{#3}}
\theoremstyle{citing}

\begin{document}

\title[]{Cup products with volume classes \\ in bounded cohomology}

\author[Michelle Bucher]{Michelle Bucher}
\address{University of Geneva, Switzerland}
\email{michelle.bucher@unige.ch}

\author[Roberto Frigerio]{Roberto Frigerio}
\address{Dipartimento di Matematica, Universit\`a di Pisa, Italy}
\email{roberto.frigerio@unipi.it}



\begin{abstract}
Let $M$ be a complete Riemannian $n$-manifold, $n\geq 2$, and suppose that either $M$ is negatively curved, or that it is a locally symmetric space of non-compact type. 
Via integration over straight simplices, every top-degree bounded differential form $\omega\in \Omega_b^n(M)$ defines a bounded cohomology class  in $H^n_b(M)$. 
We prove that every such class belongs to the annihilator of the cup product in the whole graded algebra  $H^{\geq 1}_b(M)$. In the case of closed orientable surfaces, this generalizes a result of Marasco, and, via transfer, provides a new proof of the fact that the cup product of the bounded Euler class of $SL(2,\mathbb{R})$ with itself vanishes. 
Our result also applies to degree-3 bounded classes constructed by Soma and Farre, thus showing that, for surface groups and for free groups,  the annihilator of the cup product 
is infinite dimensional also in degree 3.
\end{abstract}

\maketitle

\section{Introduction}
Let $M$ be either a complete negatively curved Riemannian $n$-manifold, or a complete locally symmetric $n$-dimensional space, $n\geq 2$.
It is known that, when $M$ has negative curvature bounded away from 0 and finite volume, 
 the bounded cohomology $H^k_b(M)$ of $M$ with real coefficients is infinite dimensional in degrees $2$ and $3$, 
but computing the bounded cohomology module $H^k_b(M)$  
for $k\geq 4$ is usually a challenging task: for example, 
it is not known if $H^4_b(M)$ vanishes or not when $M$ is a hyperbolic surface, or a hyperbolic $3$-manifold with non-elementary fundamental group.

A natural question is whether non-trivial bounded cohomology classes in degree $k\geq 4$ may be constructed by taking the cup product of lower degree classes. 
Our assumptions on $M$ ensure that $M$ supports a straightening procedure for simplices such that top-dimensional straight simplices
have uniformly bounded volumes.
Therefore, if $n=\dim M \geq 2$ and $\omega\in \Omega^n_b(M)$ 
is a top--degree bounded differential form, then
by integrating $\omega$  on straight simplices, one may construct a singular bounded cocycle which defines in turn a bounded cohomology class $[\omega]_b  \in H^n_b(M)$.
Our main theorem, which applies to complete negatively curved Riemannian manifolds with curvature bounded away from 0 and  complete locally symmetric spaces of dimension $\geq 2$,  shows that the cup product of such a class with any other bounded cohomology class always vanishes (see Section \ref{str2:sub} for the definition of a bounded straightening):
\begin{thm_intro}\label{main:thm} 
Let $M$ be a  Riemannian $n$-dimensional manifold admitting a bounded straightening. For every $k\geq 2$, every $\omega\in \Omega_b^n(M)$ and every class $\beta \in H^k_b(M)$ it holds that 
$$
[\omega]_b \cup \beta = 0 \quad \text{in}\ H^{n+k}_b(M)\ .
$$
\end{thm_intro}
In other words, bounded cohomology classes defined via integration of top--degree differential forms belong to the annihilator of the cup product (strictly speaking, non-trivial classes in $H^0_b(M)\cong \mathbb{R}$ are not annihilated by any non-trivial bounded cohomology class; here and henceforth, we will say that a class belongs to the annihilator of the cup product if it annihilates every bounded cohomology class of positive degree).

There exist by now quite a number of results on the vanishing of cup products in bounded cohomology. 
Marasco proved in~\cite{Marasco} that, for a closed negatively curved manifold of any dimension, 
bounded cohomology classes defined via integration of {degree $2$ exact} differential forms belong to the annihilator of the cup product.
When $M$ is a closed orientable surface, our Theorem~\ref{main:thm} generalizes Marasco's result by including also non-exact differential forms, like the area form.
For complete orientable hyperbolic surfaces, if $\omega$ is the area form, then $[\omega]_b$ is a nonzero multiple of the bounded Euler class, hence Theorem~\ref{main:thm} implies the following:

\begin{cor_intro}\label{euler:cor}
Let $\Sigma=\mathbb{H}^2/\Gamma$ be a complete orientable hyperbolic surface, and let $e_b^\Sigma \in H^2_b(\Sigma)\cong H^2_b(\Gamma)$ be the bounded Euler class associated to
the action of $\Gamma$ on $\partial \mathbb{H}^2$. Then
$$
e_b^\Sigma  \cup \beta = 0 \quad \text{in}\ H^{2+k}_b(\Sigma)
$$ 
for every class $\beta\in H^k_b(\Sigma)$, $k\geq 1$.
\end{cor_intro}
When $\Gamma$ is a uniform lattice in $G=\mathrm{SL}(2,\mathbb{R})$ (i.e., when $\Sigma$ is closed), the restriction map from the continuous bounded cohomology 
$H^\bullet_{cb}(\mathrm{SL}(2,\mathbb{R}))$ to $H^\bullet_b(\Gamma)$ 
admits a left-inverse  \emph{transfer map}, hence it
 is injective (see, e.g.,~\cite[Proposition 2.4.1]{BM02}). Moreover, 
via transfer, the class $e_b^\Sigma$ is taken to the bounded Euler class $e_b\in H^2_{cb}(\mathrm{SL}(2,\mathbb{R}))$. 
Therefore, Corollary~\ref{euler:cor}
provides a new proof of the following theorem,
which was proved independently by Morita~\cite{morita} and by Burger and Monod~\cite{bmcup}:

\begin{cor_intro}\label{eulerPSL:cor}
Let $e_b\in H^2_{cb}(\mathrm{SL}(2,\mathbb{R}))$ be the bounded Euler class. Then
$$
e_b\cup e_b=0 \quad \text{in}\ H^{4}_{cb}(\mathrm{SL}(2,\mathbb{R}))\ .
$$
\end{cor_intro}

Further vanishing results for the cup product of bounded classes coming from differential forms were recently obtained in~\cite{FriMil}. It is maybe worth mentioning that,
at least in degree 2, such classes span an infinite-dimensional subspace of $H^2_b(M)$ for every closed negatively curved $n$-manifold $M$, $n\geq 2$ 
(see \cite{bargeghys,Battista}). 

In the past years
 a lot of attention has been devoted to the case of the free
group, for which it has been shown that cup products of bounded cohomology classes with degree 2 classes vanish in a wide variety of cases
(see~\cite{BucherMonod,Heuer,Bucher2}). Recall that, if $X$ is any topological space, then the bounded cohomology of $X$ is isometrically isomorphic 
to the bounded cohomology of the fundamental group of $X$, via an isomorphism that preserves the cup product.
If $F_r$ is the free group of rank $r\geq 2$ and $M$ is  
any complete non-compact hyperbolic surface of finite area with $\pi_1(M)\cong F_r$, then it is shown in~\cite{DAFH} that the map
$\Omega^2_b(M)\to H^2_b(M)\cong H^2_b(F_r)$ which sends an element $\omega\in \Omega^2_b(M)$ to the class $[\omega]_b$ is injective.
Therefore, Theorem~\ref{main:thm} provides yet another infinite-dimensional subspace of $H^2_b(F_r)$ contained in the annihilator of the cup product.

Bounded differential forms have been extensively exploited also to exhibit non-trivial elements in the third bounded cohomology group 
of free groups and of surface groups. In fact, if $\Gamma$ is either a free group or a surface group, then infinite-dimensional subspaces of
$H^3_b(\Gamma)$ may be constructed by taking classes of the form $[\omega]_b$, where $\omega$ is a bounded differential $3$-form on a (geometrically infinite)
hyperbolic $3$-manifold $M$ with $\pi_1(M)\cong \Gamma$. We refer the reader e.g.~to~\cite{soma1, soma2,  farre0, farre} for more details.
Therefore, our Theorem~\ref{main:thm} also implies the following:

\begin{cor_intro}\label{soma:cor}
Let $\Gamma$ be either a surface group or a non-abelian free group. Then the subspace of the third cohomology group $H^3_b(\Gamma)$ given by elements
in the annihilator of the cup product is infinite dimensional.
\end{cor_intro}

\subsection*{Plan of the paper} In Section~\ref{bounded_cohomology:sec} we recall the definition of bounded cohomology. In Section \ref{str2:sub}, we introduce the notion of (bounded) straightening, give examples of such straightenings and  provide the precise definition of the cohomology class  $[\omega]_b$ involved in our main theorem. Section~\ref{filling:sec} is devoted to  the proof of Theorem~\ref{main:thm}.

\subsection*{Acknowledgements}
The  first author is partially supported by the SNSF and the second author is partially supported by INdAM through GNSAGA.

\subsection*{AI disclosure} In a previous version of the paper, the authors proved Theorem~\ref{main:thm} and Corollary~\ref{euler:cor} 
under the additional hypothesis that $M$ be compact, without any use of AI. This was sufficient to get also Corollary~\ref{eulerPSL:cor}. We submitted that version to ChatGPT 6.1 Sol,
asking whether our proofs were correct and/or could be simplified, and AI observed that (after a minor modification) our proof would work also when replacing compactness with 
completeness in the hypothesis of our main theorem. This allowed us to get also Corollary~\ref{soma:cor}.

\section{Bounded cohomology}\label{bounded_cohomology:sec}
Let $M$ be a topological space (in fact, we will be interested only in the case when $M$ is a manifold). 
For every $k \in \mathbb{N}$, we denote by $\Delta^k$ the standard $k$-simplex, which is the convex hull of the canonical basis of $\mathbb{R}^{k+1}$, and by
$\Simpl_k(M)$  the set of singular $k$-dimensional simplices of $M$, i.e., the set whose elements are continuous maps from $\Delta^k$ to $M$.
We denote by $\ChainSing_\bullet (M)$ the singular chain complex of $M$ with real coefficients, which is given, in degree $k$, by the real vector space
with basis $\Simpl_k(M)$.
Henceforth, we implicitly understand that (co)chains and (co)homology groups are always with real coefficients.

\subsection{Bounded cohomology}

The \emph{bounded cohomology} $H^\bullet_b(M)$ of $M$ is the cohomology of the cocomplex
\[
	0\longrightarrow \CochainSing^0(M)\longrightarrow \CochainSing^1(M)\longrightarrow \CochainSing^2(M)\longrightarrow\cdots \,,
\]
where $\CochainSing^k(M)$ denotes the space of bounded singular $k$-cochains of $M$, the differential maps are the restrictions of the coboundary maps of the usual singular cochain complex $C_\mathrm{Sing}^\bullet(M)$ and a cochain $b\in C^k_\mathrm{Sing}(M)$ is \emph{bounded} if
\[
	\lVert b \rVert_\infty =
	\sup \bigl\{ | b(\sigma) |, \; \sigma \in \Simpl_k(M)
	\bigr\} < \infty \ .
\]

\subsection{Vertex cochains}
Let $M$ be any topological space, and denote by $ \widetilde{M}$ its universal covering. If $\Gamma\cong \pi_1(M)$ is the automorphism group of the universal covering,
we define the set of \emph{vertex $k$-simplices} on  $\widetilde{M}$ simply as
$$\Comb_k(\widetilde{M})=\widetilde{M}^{k+1}\ .$$ 
We denote by $\Chain_\bullet (\widetilde{M})$ the vertex chain complex on $\widetilde{M}$ with real coefficients, which is given, in degree $k$, by the real vector space
with basis $\Comb_k(\widetilde{M})$. The differential of the complex $\Chain_\bullet (\widetilde{M})$ is given, in degree $k$ and on the basis $\Comb_k(\widetilde{M})$, by the formula
$$
\partial_k(x_0,\dots,x_{k})=\sum_{i=0}^{k} (-1)^i (x_0,\dots, \widehat{x}_i,\dots,x_{k})\ .
$$

Following~\cite[Section 2.3]{Gromov82},
we define the complex of \emph{vertex}\footnote{Gromov calls them \emph{straight} cochains, but we find this terminology confusing in our context of straightenings.} \emph{cochains} on $M$ by setting, for every $k\in\mathbb{N}$, 
$$
\Cochain^k(\tilM)=\ell^\infty (\widetilde{M}^{k+1})\ ,
$$
where $\ell^\infty (\widetilde{M}^{k+1})$ denotes the space of real-valued bounded functions on $\widetilde{M}^{k+1}$, and
$$
\Cochain^k(M)=\ell^\infty (\widetilde{M}^{k+1})^\Gamma\ ,
$$
where $\Gamma$ is understood to act diagonally on $\widetilde{M}^{k+1}$. The differential 
$$\delta^k: \Cochain^k(M)\longrightarrow \Cochain^{k+1}(M)$$
is induced by $\partial_{k+1}$. 

There is a natural $\Gamma$-equivariant map 
$$\nu_k: \ChainSing_k(\widetilde{M})\longrightarrow \Chain_k(\widetilde{M})$$
defined on the basis of singular simplices $\sigma\in  \Simpl_k(\widetilde{M})$ by sending $\sigma$ to its $k+1$ ordered vertices. (This provides a justification for our choice of terminology, as we think of a vertex chain in $\Chain_k(\widetilde{M})$ as a linear combination of simplices which only remember the positions of their vertices.) Given the obvious identification between $\CochainSing^k(M)$ and $\CochainSing^k (\widetilde{M})^\Gamma$, the $\Gamma$-equivariant chain map $\nu_\bullet$ induces a map of cocomplexes 
$$
\nu^\bullet\colon \Cochain^\bullet(M)\longrightarrow \CochainSing^\bullet (\widetilde{M})^\Gamma\cong \CochainSing^\bullet(M)\ .
$$

Originally, Gromov required vertex cochains to be Borel measurable (or even continuous) bounded maps from $\widetilde{M}^{n+1}$ to $\mathbb{R}$. It was then proved  by L\"oh in~\cite[Theorem 3.5]{Loehprop} that continuous vertex cochains compute the bounded cohomology of $M$, i.e., that
the restriction to continuous vertex cochains of the cochain map $\nu^\bullet$ induces an (isometric) isomorphism in bounded cohomology. 

To have an explicit chain homotopy inverse of $\nu_\bullet$,  it is necessary to have a procedure that assigns to any $(k+1)$-tuple of points in $\widetilde{M}$ a singular simplex in $\widetilde{M}$ (in a $\Gamma$-equivariant way). This is precisely the purpose of a straightening and is discussed in the next section. For manifolds satisfying the very restrictive hypothesis of supporting a straightening, we will recover a well known variant of L\"oh's isomorphism, namely that $\nu^\bullet$ induces an isomorphism on bounded cohomology with no regularity requirement (except boundedness) for the maps from $\widetilde{M}^{n+1}$ to $\mathbb{R}$.

\section{Straightenings}\label{str2:sub}
Let $M$ be a smooth manifold. As before let $\widetilde{M}$ be its universal covering and $\Gamma\cong \pi_1(M)$ its fundamental group.

\begin{defn}\label{defn:straightening}
	
	A \emph{straightening} $\STR_\bullet$ on $M$ is a sequence of $\Gamma$-equivariant maps
	\[ \STR_k \colon \Comb_k(\widetilde{M}) \to \Simpl_k(\widetilde{M}), \]
	where $k$ varies in $\mathbb{N}$, with the following properties:
	\begin{enumerate}
		\item For every $(v_0,\dots,v_k)\in \Comb_k(\tilM)$, its straightening $\STR_k(v_0,\dots,v_k)\colon \Delta^k \to \tilM$ is of class $C^1$;
		\item The $i$-th face of the straightening   $\STR_k(v_0,\dots,v_k)$ of any $(v_0,\dots,v_k)\in \Comb_k(\tilM)$ is equal to the straightening $\STR_{k-1}(v_0,\dots,\widehat{v_i},\dots,v_k)$ of the $i$-th face $(v_0,\dots,\widehat{v_i},\dots,v_k)\in \Comb_{k-1}(\tilM)$ of $(v_0,\dots,v_k)$;
		\item For any $\sigma\in \Simpl_k(\widetilde{M})$, there is a family of homotopies $h_\sigma \colon \Delta^k \times [0,1] \to \widetilde{M}$ connecting $\sigma$ and $\STR_k(\nu_k(\sigma))$, where $k$ varies in $\mathbb{N}$; this family of homotopies is $\Gamma$-equivariant and compatible with restriction to faces, i.e., the homotopy associated to a  face of $\sigma$ 
		is equal to the restriction of $h_\sigma$ to the same face times $[0,1]$.
		\end{enumerate}
If furthermore the straightening satisfies the following additional property, then it is called a \emph{bounded straightening}:
		\begin{enumerate}
		\item[(4)] Let $n=\dim M=\dim \tilM$. There exists a constant $D>0$
such that, for every \emph{straight} simplex  $s=\STR_n(v_0,\dots,v_n)\colon \Delta^n\to \tilM$,
$$
\int_{\Delta^n} |\det ds_p|\, d\mu(p) \leq  D\ ,
$$
where $ds_p$ denotes the differential of $s$ at $p$, 
the determinant of $ds_p$ is computed by taking orthonormal
frames in $T_p \Delta^n$ and in $T_{s(p)} \tilM$, and
$d\mu$ is the usual Euclidean measure on $\Delta^n$.
 
	\end{enumerate}
\end{defn}

The linear extension of a straightening provides a map 
$$\STR_\bullet\colon \Chain_\bullet(\tilM) \to \ChainSing_\bullet(\widetilde{M}),$$ 
which is a $\Gamma$-equivariant chain map by Property (2). It is easy to see that 
$$ \nu_\bullet\circ \STR_\bullet:\Chain_\bullet(\tilM) \longrightarrow \Chain_\bullet(\tilM)  $$
is $\Gamma$-equivariantly chain homotopic to the identity. (In the many cases where $\STR_0$ is the identity, it is already the identity. In the general case an explicit homotopy is easy to exhibit.)  Furthermore by Property (3) we also have that 
$$\STR_\bullet \circ \nu_\bullet: \ChainSing_\bullet(\widetilde{M})\longrightarrow \ChainSing_\bullet(\widetilde{M})$$
is $\Gamma$-equivariantly chain homotopic to the identity, see e.g.~\cite[Proposition 8.11]{Frigeriobook}. Both maps are norm non increasing for the $\ell^1$-norms at chain level and thus induce maps on the corresponding cochain complexes, which are mutually inverses at cohomology level: 

\begin{lemma}\label{straight:lem}
Suppose that the manifold $M$ supports a straightening, according to Definition~\ref{defn:straightening}. Then the cochain maps 
$$\nu^\bullet\colon \Cochain^\bullet(M) \longrightarrow  \CochainSing^\bullet({M}) $$ 
and 
$$\STR^\bullet\colon  \CochainSing^\bullet({M}) \longrightarrow \Cochain^\bullet(M)$$ 
induce isometric isomorphisms of graded algebras on bounded cohomology. In particular the bounded cohomology of $M$ can be computed as the cohomology of the vertex cochain cocomplex. 
\end{lemma}

When considering the bounded cohomology $H^\bullet_b(M)$ of $M$ we will henceforth work with the vertex cocomplex. 

Observe that, in contrast to our  Main Theorem \ref{main:thm}, for this lemma to hold it is not necessary to have a bounded straightening. 

We finish this section by giving several examples of (bounded) straightenings before seeing how the existence of a bounded straightening (and here boundedness will play a crucial role)  is used to ensure the existence of a bounded representative of the top degree volume class (Subsection \ref{subsect bdd vol}).

\subsection{The geodesic straightening}\label{geodesic:subsec}
Suppose that $\dim M\geq 2$, and that $M$ is negatively curved with curvature bounded away from 0.
Since $M$ is  complete, the universal covering $\widetilde{M}$
is itself complete. Being also negatively curved,  it is uniquely geodesic, and this allows us to define a straightening inductively on $\widetilde{M}$ by coning a $(k-1)$-singular simplex with vertices $v_1,\dots,v_k$ with respect to the geodesic segments from $v_0$.

Property~(1) in Definition~\ref{defn:straightening}, i.e., the fact that the straight simplex 
is indeed of class $C^1$ is in fact false in general, but can easily be achieved by reparametrising the simplices coherently with respect to the faces.   Property (2) readily follows from the definition, while
property (3) is shown e.g.~in~\cite[Proposition 8.11]{Frigeriobook}. The crucial property (4), i.e., the fact that 
 the geodesic straightening has uniformly bounded volumes, is proved in~\cite{inoue_gromov_1982}. 
 
\subsection{The barycentric straightenings}\label{bary:sub} Suppose now that $\dim M \geq 3$, and that $M$ is a locally symmetric space of non-compact type. 
Given $(x_0,\dots,x_k)\in \widetilde{M}^{k+1},$ a straightening can be defined by the map
$$\begin{array}{rcl}
\Delta^k&\longrightarrow& \widetilde{M}\\
(t_0,\dots,t_k)&\longmapsto & \mathrm{the \ unique \ minimizer \ of \ }  \widetilde{M}\ni x\mapsto \sum_{i=0}^k t_i \, d(x,x_i)^2.
\end{array}$$ 
The fact that the minimizer is unique follows from the convexity properties of the distance. It is  straightforward to check that this defines a straightening, but its boundedness, while likely, seems hard to establish. The unique minimizer in this construction can also be seen as the Fr\'echet barycenter  of the measure $\sum_{i=0}^k t_i \delta_{x_i}$, where $\delta_{x_i}$ denotes the Dirac measure at $x_i\in \tilM$. Replacing the Dirac measures of the vertices by Patterson-Sullivan measures on the boundary at infinity of the symmetric spaces, and the Fr\'echet barycenter by the Busemann barycenter, Besson, Courtois and Gallot~\cite{BCG} and, for higher rank symmetric spaces, Connell and Farb~\cite{CF}, considered a barycentric straightening which is a straightening according to our Definition~\ref{defn:straightening}. Lafont and Schmidt used a slight variant of this straightening in~\cite{LFS} for which they could further show boundedness (our Property (4)) except in the case of locally symmetric spaces with local $\mathbb{H}^2$ or $\mathrm{SL}(3, \mathbb{R})/\mathrm{SO}(3)$ factors, for which one needs
to appeal, respectively, to the geodesic straightening described above (which, for hyperbolic manifolds, coincides up to reparametrization with the barycentric straightening), or to the projective straightening described in the next example. 

A bounded barycentric straightening may be defined and satisfies all the properties described in Definition~\ref{defn:straightening} also
when $M$ is a closed  $n$-dimensional negatively curved manifold, $n\geq 3$ (see~\cite{connell_positivity_2019} for a proof). 


\subsection{The projective straightening} If $M$ is a convex projective manifold, then its universal cover is a properly convex domain in some affine chart of the real projective $n$-space $\mathbb{P}^n(\mathbb{R})$. It is not uniquely geodesic in general, but one can nevertheless mimic the inductive geodesic straightening for the favorite geodesics which are given by (reparametrizations of) affine segments in the affine chart \cite[Section 8]{BFRProjective}. This is always a straightening. It is known to be a bounded straightening for manifolds that are locally isomorphic to $\mathrm{SL}(3, \mathbb{R})/\mathrm{SO}(3)$ \cite{BucherSL}. These straightenings are exploited in  \cite{BFRProjective} to give the inequality $|\chi(M)|\leq ((n+2)/4) \|M\|$  between the Euler characteristic and the simplicial volume of any closed, oriented, convex projective $n$-dimensional manifold $M$.

Note that in the case of hyperbolic manifolds, the geodesic, the barycentric and the projective straightening coincide, up to reparametrization.

\begin{rmk}
By the very definition, the map $\str_0$ obtained via the geodesic or the projective straightening  is the identity. In other words, 
the geodesic straightening of $(v_0,\dots,v_k)$ has vertices $(v_0,\dots,v_k)$. This is not true in general
for the barycentric straightening, and this is the reason why, 
in \Cref{defn:straightening}, we do  {not} require $\str_0$ to be the identity. 
\end{rmk}

\subsection{Bounded cohomology classes from differential forms}\label{subsect bdd vol} Let as before $M=\tilM/\Gamma$ be a complete $n$-dimensional manifold. 
For every smooth differential $n$-form $\omega\in \Omega^n(M)$ we set
$$
\|\omega\|_\infty=\sup \{|\omega_p(e_1,\dots, e_n)|\, ,\ p\in M\, ,\ (e_1,\dots,e_n)\ \text{orthonormal}\ \text{frame in}\ T_p M\}\ .
$$
We say that $\omega$ is \emph{bounded} if $\|\omega\|_\infty<\infty$ (this condition always holds if $M$ is compact), and we denote by
$\Omega^n_b(M)\subseteq \Omega^n(M)$ the subspace of bounded $n$-forms. If $M$ is oriented with volume form $d\text{Vol}$, then
any $n$-form $\omega\in \Omega^n(M)$ is of the form $\omega=f\, d\text{Vol}$ for some smooth function $f\colon M\to \mathbb{R}$,
and $\omega$ is bounded exactly when $f$ is bounded. 

Let $\omega \in \Omega_b^n(M)$ be a bounded differential form. Being top dimensional, $\omega$ is automatically closed. Even if one could disregard the fact that singular simplices $\sigma\in \ChainSing_n(M)$ are not smooth in general, the singular cocycle defined by 
$$\sigma \longmapsto \int_\sigma \omega$$ 
is not bounded unless $\omega=0$. However, in the presence of a bounded straightening 
\[ \STR_\bullet \colon \Comb_\bullet(\tilM) \longrightarrow \Simpl_\bullet(\widetilde{M}) \]
on $M$, we can instead consider the singular cocycle defined by 
\begin{equation}\label{equ: sing cocycle} 
\sigma \longmapsto \int_{p\circ \str_n (\nu_{n}(\widetilde{\sigma}))} \omega\ ,
\end{equation}
where $\widetilde{\sigma}$ is any lift of $\sigma$ to $\widetilde{M}$ and $p:\widetilde{M}\rightarrow M$ denotes the universal covering projection. Property (4) of the definition \ref{defn:straightening} of bounded straightenings, together with the boundedness of $\omega$, ensures that this is indeed bounded. Likewise, letting $\widetilde{\omega}\in\Omega^n(\widetilde{M})$ be the pull-back of $\omega$ to $\widetilde{M}$, we can consider the corresponding vertex bounded cocycle 
$$\begin{array}{rccl}
I(\widetilde{\omega}):&\widetilde{M}^{n+1}&\longrightarrow & \mathbb{R}\\
&(x_0,\dots,x_n)&\longmapsto &  \int_{ \str_n (x_0,\dots,x_n)} \widetilde{\omega} \ . 
\end{array}$$
Note that the claim that these two cochains are indeed cocycles follows from the combination of Stokes' Theorem and the facts that the forms $\omega$ and $\widetilde{\omega}$ are closed,  and that $\str_\bullet$ is a chain map.

Taking the vertex cochain complex $\Cochain^\bullet(M)=\ell^{\infty}(\widetilde{M}^{\bullet +1})^\Gamma$ for computing the bounded cohomology $H^\bullet_b(M)$, we can thus define a map 
$$\begin{array}{rcl}
 \Omega_b^n(M)&\longrightarrow & H^n_b(M)\\
 \omega&\longmapsto & [I(\widetilde{\omega})]\ .
\end{array}$$
Note that we could  likewise have defined this map in the singular cochain complex using the cocycle (\ref{equ: sing cocycle}). This is the cohomology class which we denoted by 
$[\omega]_b$ in the introduction. Finally observe that the latter map may depend on the bounded straightening. 

The procedure described in this section is the original motivation behind the introduction of bounded straightenings as it  shows that, for closed orientable manifolds admitting a bounded straightening, the top dimensional volume form in the singular cohomology $H^n_\mathrm{Sing}(M)$ (which is isomorphic to the de Rham cohomology) can be represented by a bounded cocycle. This  is equivalent to the strict positivity of the simplicial volume of $M$. 

\section{Cup product of top--degree forms}\label{filling:sec}
Let us now come back to the general situation when $M=\widetilde{M}/\Gamma$ is  a complete Riemannian $n$-manifold, where $\Gamma<\text{Isom}(\widetilde{M})$ is isomorphic to $\pi_1(M)$, and
let $p\colon \widetilde{M}\to M$ be the universal covering projection. 
Throughout this section, we assume that $M$ supports a bounded straightening 
$$\str_\bullet\colon \Comb_\bullet(\tilM) \longrightarrow \Simpl_\bullet(\widetilde{M})\ ,$$
and we fix a bounded differential form $\omega\in \Omega^n_b(M)$. 
We set as usual $\widetilde{\omega}=p^*(\omega)\in \Omega^n_b(\tilM)$.
Our goal is to prove that 
the class $[\omega]_b\in H^n_b(M)$ represented by the vertex cocycle $I(\widetilde{\omega})\in \Cochain^n(M)$ belongs to the annihilator
of the cup product.

We first express $I(\widetilde{\omega})$ 
as an infinite sum of (non $\Gamma$-invariant) cocycles 
$I_\gamma(\widetilde{\omega})\in \Cochain^n(\tilM)$
indexed by the elements $\gamma$ of $\Gamma$.
Let $F\subseteq \tilM$ be a fundamental domain for the action of $\Gamma$ on $\tilM$.
For any $\gamma\in\Gamma$ and every $(x_0,\dots,x_n)\in\widetilde{M}^{n+1}$ we set 
$$
I_\gamma (\widetilde{\omega})(x_0,\dots,x_n)=\int_{\str(x_0,\dots,x_n)} \chi_{\gamma}\, \widetilde{\omega} \in \mathbb{R}\ ,
$$
where $\chi_{\gamma}\colon \widetilde{M}\to\mathbb{R}$ is the characteristic function of $\gamma F$.
Since the family $\{\gamma F\}_{\gamma\in\Gamma}$ is locally finite and the support of
$\str(x_0,\dots,x_n)$ is compact, for every fixed $(x_0,\dots,x_n)\in\widetilde{M}^{n+1}$
we have  $I_\gamma (\widetilde{\omega})(x_0,\dots,x_n)\neq 0$ only for a finite number of elements $\gamma\in\Gamma$. Therefore,  the infinite sum $\sum_{\gamma \in\Gamma} I_\gamma (\widetilde{\omega})$ 
provides a well defined map $\widetilde{M}^{n+1}\to \mathbb{R}$. By definition of $F$ we  have $\sum_{\gamma \in \Gamma} \chi_\gamma\widetilde{\omega}=\widetilde{\omega}$
almost everywhere on $\widetilde{M}$, hence 
$$
\sum_{\gamma \in \Gamma} I_\gamma (\widetilde{\omega}) = I(\widetilde{\omega})\ .
$$

\begin{lemma}\label{Iomega}
\begin{enumerate}
\item For every $\gamma\in\Gamma$, the vertex cochain $I_\gamma(\widetilde{\omega})$ is a cocycle.
\item  For every $\gamma_0,\gamma\in\Gamma$, we have $\gamma_0\cdot I_\gamma (\widetilde{\omega})=I_{\gamma_0 \gamma}(\widetilde{\omega})$.
\item We have
$$\sup_{(x_0,\dots,x_n)\in\widetilde{M}^{n+1}} \sum_{\gamma \in \Gamma} |I_\gamma(\widetilde{\omega})(x_0,\dots, x_n)| \leq D\|\omega\|_\infty\ ,$$
where $D$ is the constant appearing in the definition of straightening with bounded volumes (in particular,
$\|I_\gamma(\widetilde{\omega})\|_\infty\leq D\|\omega\|_\infty$ for every $\gamma\in\Gamma$, hence $I_\gamma(\widetilde{\omega})$ is indeed bounded).
\end{enumerate}
\end{lemma}

\begin{proof}
(1): For every $\gamma\in\Gamma$, the Borel measurable form $\chi_\gamma\, \widetilde{\omega}$, being top--dimensional,
vanishes on the boundary of any $(n+1)$-singular chain, 
hence 
\begin{align*}
0 & = \int_{\partial \str(x_0,\dots,x_{n+1})}\chi_\gamma\, \widetilde{\omega} \ = \
\sum_{i=0}^{n+1} (-1)^i  \int_{\str(x_0,\dots,\widehat{x}_i,\dots, x_{n+1})}\chi_\gamma\, \widetilde{\omega}\\ & =
\sum_{i=0}^{n+1} (-1)^i I_\gamma (\widetilde{\omega})(x_0,\dots,\widehat{x}_i,\dots, x_{n+1}) \ = \
\delta I_\gamma (\widetilde{\omega}) (x_0,\dots,x_{n+1})\ .
\end{align*}

(2): 
Using that the differential form $\widetilde{\omega}$ is $\Gamma$-invariant, for every $(x_0,\dots,x_n)\in \widetilde{M}^{n+1}$ we get
\begin{align*}
(\gamma_0\cdot I_{\gamma})(\widetilde{\omega})(x_0,\dots, x_n)& =
I_{\gamma}(\widetilde{\omega})(\gamma_0^{-1}x_0,\dots,\gamma_0^{-1} x_n) =\int_{\str(\gamma_0^{-1}x_0,\dots,\gamma_0^{-1}x_n)} \chi_{\gamma}\, \widetilde{\omega} \\ & =
\int_{\str(x_0,\dots,x_n)} \chi_{\gamma_0\gamma}\, \widetilde{\omega} 
\ = \ I_{\gamma_0\gamma}(\widetilde{\omega})(x_0,\dots,x_n)\ .
\end{align*}

(3): Let $\widetilde{\omega}=f\, d\vol$, where $d\vol$ denotes the volume differential form on $\widetilde{M}$ and
$f\colon \widetilde{M}\to \mathbb{R}$ is a smooth function (even when $M$ is not orientable, we can fix an auxiliary orientation on
the simply connected manifold $\tilM$). By definition of $\ell^\infty$-norm of a differential form 
(and since $\|\widetilde{\omega}\|_\infty=\|\omega\|_\infty$) we
have $|f(x)|\leq \|\omega\|_\infty$ for every $x\in\tilM$. 
Let $(x_0,\dots,x_n)\in\widetilde{M}^{n+1}$ be fixed, and let for ease of notation $s=\str(x_0,\dots,x_n)$.
For every $\gamma\in\Gamma$ we set
$B_\gamma=s^{-1}(\gamma F)\subseteq \Delta^n$. Then
\begin{align*}
 |I_\gamma (\widetilde{\omega})(x_0,\dots,x_n)| &=\left|\int_{\str(x_0,\dots,x_n)} \chi_{\gamma}\, \widetilde{\omega}\right|=\left| \int_{B_\gamma}  f(s(p)) \det ds_p\, d\mu(p) \right|\\
 & \leq \int_{B_\gamma}  |f(s(p))|\cdot  |\det ds_p|\, d\mu(p)\leq 	\|\omega\|_\infty \int_{B_\gamma}   |\det ds_p|\, d\mu (p)\ .
\end{align*}
Now, by construction $\Delta^n=\bigcup_{\gamma\in\Gamma} B_\gamma$, and
for every pair of distinct elements $\gamma\neq \gamma'$ of $\Gamma$ the set $B_\gamma\cap B_{\gamma'}$
is sent by $s=\str(x_0,\dots,x_n)$ to a measure-zero subset of $\widetilde{M}$, hence
\begin{align*}
\sum_{\gamma \in \Gamma} |I_\gamma(\widetilde{\omega})(x_0,\dots, x_n)|  & \leq \|\omega\|_\infty \cdot \sum_{\gamma \in\Gamma} \int_{B_\gamma}   |\det ds_p|\, d\mu (p)\\ &= \|\omega\|_\infty \cdot\int_{\Delta^n} |\det ds_p|\, d\mu (p)\leq  D \|\omega\|_\infty \ .
\end{align*}
\end{proof}

 We are now ready to show that
$$
[\omega]_b \cup \beta= 0 \quad \text{in}\ H^{n+k}_b(M)\ ,
$$ 
\text{for every} $\beta\in H^k_b(M)$,
thus concluding the proof of  Theorem~\ref{main:thm}.
  We use Lemma \ref{straight:lem} to identify the bounded cohomology $H^\bullet_b(M)$ of $M$ with the cohomology of the vertex cochain complex. 
  In this complex, the class $[\omega]_b$ is represented by the cocycle $I(\widetilde{\omega})$, and a generic
  element
  $\beta \in H^k_b(M)$ is represented by a bounded vertex cocycle
$b \in \ell^\infty(\widetilde{M}^{k+1})^\Gamma$. In order to conclude  we need to exhibit a bounded $\Gamma$-invariant primitive of $I(\widetilde{\omega})\cup b$ in $\ell^\infty(\widetilde{M}^{n+k})^\Gamma$.

 To this aim, for every $\gamma\in\Gamma$ let us first define the cone $c_\gamma (b)\in \ell^\infty (\widetilde{M}^{k})$ of $b$ as follows: we fix a basepoint $o\in \widetilde{M}$, and we set
$$
c_\gamma(b)(x_1,\dots,x_k)=b(\gamma o,x_1,\dots, x_k)\quad \text{for every}\ (x_1,\dots,x_k)\in \widetilde{M}^k\ .
$$
By construction we have $\|c_\gamma(b)\|_\infty\leq \|b\|_\infty$ for every $\gamma\in\Gamma$. 
The cochain $c_\gamma(b)$ is not invariant: indeed, since $b$ is $\Gamma$-invariant, for every $\gamma_0,\gamma\in \Gamma$ we have
\begin{align*}
(\gamma_0\cdot c_\gamma(b))(x_1,\dots,x_k)& =
c_\gamma(b)(\gamma_0^{-1}x_1,\dots,\gamma_0^{-1}x_k)=b(\gamma o,\gamma_0^{-1}x_1,\dots, \gamma_0^{-1}x_k)\\ &=
b(\gamma_0\gamma  o,x_1,\dots, x_k)=c_{\gamma_0\gamma}(b)(x_1,\dots,x_k)\ ,
\end{align*}
hence $\gamma_0\cdot c_\gamma(b)=c_{\gamma_0\gamma}(b)$. 
Moreover, using that $c_\gamma (b)$ is obtained from $b$ by coning over $\gamma o$, it is immediate to check that
$\delta c_\gamma (b)=b$.

Let us set
$$
P=\sum_{\gamma \in \Gamma} I_\gamma(\widetilde{\omega})\cup c_\gamma(b)\ \in \ \ell^\infty (\widetilde{M}^{n+k})\ .
$$
We are now going to show that $P$ is (up to a sign) a bounded $\Gamma$-invariant primitive of $I(\widetilde{\omega})\cup b$, thus
concluding the proof of Theorem~\ref{main:thm}.

We first show that $P$ is bounded. Recall that, for any fixed $(n+1)$-tuple $(x_0,\dots,x_n)\in\widetilde{M}^{n+1}$, we have
$I_\gamma(\widetilde{\omega})(x_0,\dots,x_n)\neq 0$ only for a finite number of elements $\gamma\in\Gamma$, hence 
for every $(x_0,\dots,x_{n+k-1})\in\widetilde{M}^{n+k}$ the sum defining 
$P(x_0,\dots,x_{n+k-1})$ is finite. 
Thanks to~Lemma~\ref{Iomega}--(3), we have
\begin{align*}
|P(x_0,\dots,x_{n+k-1})| & \leq \sum_{\gamma \in \Gamma} | I_\gamma(\widetilde{\omega})(x_0,\dots, x_n)|\cdot |c_\gamma(b) (x_n,\dots, x_{n+k-1})|
\\ & \leq \|b\|_\infty \sum_{\gamma \in \Gamma} | I_\gamma(\widetilde{\omega})(x_0,\dots, x_n)| \leq D\|b\|_\infty \|\omega\|_\infty\ ,
\end{align*}
hence $P$ is indeed bounded. 

In order to show that $P$ is $\Gamma$-invariant it suffices to observe that, if $\gamma_0\in\Gamma$, then
$$
\gamma_0\cdot P=\sum_{\gamma \in \Gamma} (\gamma_0\cdot I_\gamma(\widetilde{\omega})) \cup (\gamma_0\cdot c_\gamma (b))=
\sum_{\gamma \in \Gamma} I_{\gamma_0\gamma}(\widetilde{\omega}) \cup c_{\gamma_0\gamma} (b)=P\ ,
$$
where in the second equality we have exploited Lemma~\ref{Iomega}--(2), and the fact that $\gamma_0\cdot c_\gamma (b)=c_{\gamma_0\gamma} (b)$.

Finally, recall that for every $\gamma\in\Gamma$ we have $\delta c_{\gamma}(b)=b$, while  $\delta I_\gamma(\widetilde{\omega})=0$ 
by Lemma~\ref{Iomega}--(1), hence
\begin{align*}
\delta P& = \delta \left(\sum_{\gamma \in \Gamma} I_\gamma(\widetilde{\omega})\cup c_\gamma(b)\right)=
\sum_{\gamma \in \Gamma} \left( \delta I_\gamma(\widetilde{\omega})\cup c_\gamma(b)\right) + (-1)^n\sum_{\gamma \in \Gamma}  \left( I_\gamma(\widetilde{\omega})\cup \delta c_\gamma(b)\right)
\\ & = (-1)^n \sum_{\gamma \in \Gamma} ( I_\gamma(\widetilde{\omega})\cup b) =  
(-1)^n \left(\sum_{\gamma \in \Gamma}  I_\gamma(\widetilde{\omega})\right) \cup b= (-1)^n I(\widetilde{\omega})\cup b\ .
\end{align*}

This concludes the proof of Theorem~\ref{main:thm}.

\bibliography{math_papers}
\bibliographystyle{alpha}
\end{document}